\documentclass[11pt]{article} 

\usepackage[utf8]{inputenc} 

\usepackage{geometry} 
\usepackage{graphicx} 

\usepackage{booktabs} 
\usepackage{array} 
\usepackage{paralist} 
\usepackage{verbatim} 
\usepackage{subfig} 
\usepackage{amsmath}
\usepackage{amssymb}
\usepackage{amsthm}
\usepackage[all]{xy}
\usepackage{leftidx}
\usepackage{amsmath,amssymb,amsfonts}
\usepackage{amsmath,amscd}
\usepackage{textcomp}

\usepackage{fancyhdr} 
\usepackage{sectsty}
\allsectionsfont{\sffamily\mdseries\upshape} 

\usepackage[nottoc,notlof,notlot]{tocbibind} 
\usepackage[titles,subfigure]{tocloft} 

\theoremstyle{plain}
\newtheorem{thm}{Theorem}

\newtheorem{defn}{Definition}
\newtheorem{remark}[defn]{Remark}

\DeclareMathOperator{\DIV}{Div}

\newcommand{\Addresses}{{
  \bigskip
  \footnotesize

  Junyan CAO, \textsc{Institut De Mathématiques De Jussieu, Université Pierre et Marie Curie,
    4, Place Jussieu, 75252 Paris Cedex 05}\par\nopagebreak
  \textit{E-mail address}, \texttt{junyan.cao@imj-prg.fr}

}}

\title{Additivity of the approximation functional of currents induced by Bergman kernels}
\author{Junyan CAO}
\date{\vspace{-5ex}}

\begin{document} 

\maketitle

\begin{abstract}
\textbf{Titre:} Additivité de la fonctionnelle d'approximation des courants induite par les noyaux de Bergman
\end{abstract}

\begin{abstract}
\textbf{Résumé:} 
Dans cette note, nous apportons une réponse positive à une question soulevée par Jean-Pierre Demailly en 2013, et démontrons l'additivité de la fonctionnelle d'approximation des courants positifs fermés de type $(1,1)$ induite par les noyaux de Bergman.
\end{abstract}

\begin{abstract}
\textbf{Abstract:}  
In this note, we give a positive answer to a question raised by Jean Pierre Demailly in 2013, and show the additivity of the approximation functional of closed positive $(1,1)$-currents induced by Bergman kernels.
\end{abstract}

\section{Introduction}

Let $X$ be a compact complex $n$-dimensional manifold. An important positive cone in complex analytic geometry is the pseudoeffective classes $\mathcal{E} (X)$, namely the subset of cohomology classes $H^{1,1} (X)$ containing a closed positive $(1,1)$-current $T= \alpha + dd^c\varphi$, where $\alpha$ is a smooth $(1,1)$-form and $\varphi$ is a quasi-psh function on $X$. 
In various geometric problems (for example, the Nadel vanishing theorem), we need to keep the information on the singularities. 
To preserve the information about the asymptotic multiplier ideal sheaves $\mathcal{I}(m\varphi)$, Demailly constructed a new cone by using in an essential way 
a Bergman kernel approximation. Before explaining this new construction, we first recall some elementary notions about quasi-psh functions.

\begin{defn}
Let $\varphi_1$, $\varphi_2$ be two quasi-psh functions on $X$. 

(1)  We say that $\varphi_1$ has analytic singularities, if locally one can write as:
$$\varphi _1= c \ln \sum_i |g_i|^2  + O (1)$$
where the $g_i$ are holomorphic functions and $c$ is a positive constant.

(2)  We say that $\varphi_1$ has less singularities than $\varphi_2$, and write as $\varphi_1 \preccurlyeq \varphi_2$, if  we have
$\varphi_1 \ge \varphi_2 +C$
for some constant $C$. 

(3) We say that $\varphi_1$ and $\varphi_2$ have equivalent singularities, and write $\varphi_1 \sim \varphi_2$, when we have both  $\varphi_1 \preccurlyeq \varphi_2$ and  $\varphi_2 \preccurlyeq \varphi_1$.
\end{defn}

We now recall briefly the constructions in \cite[Section 3]{Dem14} and recommend the reader to see \cite{Dem14} for its applications.
Let $\mathcal{S} (X)$ be the set of singularity equivalence classes of closed positive $(1,1)$-currents. 
It is naturally equipped with a cone structure.
Continuing his work \cite{Dem92} of the early 1990's on the approximation
theorem, Demailly recently defined in \cite{Dem14} another cone which has a more algebraic appearance:

\begin{defn}

For each class $\alpha \in \mathcal{E} (X)$, we define $\widehat{\mathcal{S}}_{\alpha} (X)$ as a set of equivalence classes of sequences of quasi-positive currents 
$T_k =\alpha +dd^c \psi_k$ (we suppose from now on that $\alpha$ is a smooth $(1,1)$-form on $X$) such that:

(a) $T_k =\alpha + dd^c \varphi_k \geq -\epsilon_k \omega$ with $\lim\limits_{k\rightarrow +\infty} \epsilon_k =0$.

(b) The functions $\varphi_k$ have analytic singularities and $\varphi_k \preccurlyeq \varphi_{k+1}$ for all $k$. We say that $(T_k)  \preccurlyeq_W (T_k ') $, if for every $\epsilon >0$ and $k$,
there exists $l$ such that $(1-\epsilon) T_k \preccurlyeq T_l '$. Finally, we write $(T_k) \sim (T_k ')$ when we have both $(T_k)  \preccurlyeq_W (T_k ')  $ and $(T_k ')  \preccurlyeq_W (T_k)  $,
and define $\widehat{\mathcal{S}}_{\alpha} (X)$ to be the quotient space by this equivalence relation.

(c) We set $\widehat{\mathcal{S}} (X):= \bigcup_{\alpha\in\mathcal{E} (X) } \widehat{\mathcal{S}}_{\alpha} (X)$.
\end{defn}

Let $\varphi$ be a quasi-psh function on $X$, and $(\varphi_k)$ be an Bergman kernel type approximation of $\varphi$, i.e.,
$\varphi_k \sim \frac{1}{k}\ln (\sum\limits_j |g_j|^2)$ on $U_i$, where $\{ U_i\}$ is a Stein cover of $X$ 
and $\{g_j\}$ is an orthonormal basis of $H^0 (V_i, \mathcal{O}_{V_i})$ for some Stein open set $V_i \Supset U_i$ with respect to the $L^2$ norm $\int_{V_i} |\cdot|^2 e^{-2k\varphi}$.
Let $\alpha + dd^c \varphi$ be a positive current. 
By using the comparison theorem (cf. for example \cite[Thm 2.2.1, step 3]{DPS01}), Demailly \cite[(3.1.10)]{Dem14} proved that the following map is well defined:
$$\textbf{B}:  \mathcal{S} (X)\rightarrow \widehat{\mathcal{S}} (X).$$
$$\alpha + dd^c \varphi \rightarrow (\alpha + dd^c \varphi_{2^k}) .$$
It is called here the Bergman kernel approximation functional.
 
\begin{remark}
Although we will not use it here, we should mention the following important property of the map $\textbf{B}$ (cf. \cite[Thm 2.2.1]{DPS01} or \cite[Cor 1.12]{Dem14}) : 
For every pair of positive numbers $\lambda ' > \lambda > 0$, there exists an integer $k_0 (\lambda, \lambda ') \in \mathbb{N}$ such that
$$\mathcal{I}(\lambda '  \varphi_k) \subset \mathcal{I}(\lambda   \varphi) \qquad\text{for } k\geq k_0 (\lambda, \lambda '). $$
\end{remark}

Evidently, both $\mathcal{S} (X)$ and $\widehat{\mathcal{S}} (X)$ admit an additive structure. 
\cite[Section 3]{Dem14} asked whether $\textbf{B}$ is a morphism for addition. 
In this short note, we will give a positive answer to this question. More precisely, we have

\begin{thm}[Main theorem]\label{mainthm}
Let $T_1 =\alpha_1 + dd^c\varphi$, $T_2 =\alpha_2 +dd^c\psi$ be two elements in  $ \mathcal{S} (X)$.
The we have $\textbf{B} (T_1 + T_2) =\textbf{B} (T_1 ) +\textbf{B} (T_2) $.
\end{thm}

\section{Proof of Main theorem}

\begin{proof}

In the setting of Theorem \ref{mainthm},  let $\tau_k$ (respectively $\varphi_k$, $\psi_k$)  be a Bergman kernel type approximation of $\varphi +\psi$ 
(respectively $\varphi$, $\psi$) . 
By the subadditive property of ideal sheaves $\mathcal{I}(k\varphi +k \psi) \subset \mathcal{I}(k\varphi) \mathcal{I}(k \psi) $ (\cite[Thm 2.6]{DEL00}), 
we have $\varphi_k + \psi_k  \preccurlyeq \tau_k $.
By Definition 2, to prove our main theorem, it is sufficient to prove that for every $k\in\mathbb{N}$ fixed, there exists a positive sequence $\lim\limits_{p\rightarrow +\infty}\epsilon_p = 0$, such that 
\begin{equation}\label{equivalproof}
(1-\epsilon_p) \tau_k \preccurlyeq \varphi_p + \psi_p \qquad\text{for  every }p \gg 1.
\end{equation}

For every $k\in\mathbb{N}$ fixed, there exists a bimeromorphic map $\pi: \widetilde{X}\rightarrow X$, such that 
\begin{equation}\label{add2}
\tau_k \circ \pi = \sum_i c_i \ln |s_i| + C^{\infty} \qquad\text{for some }c_i >0 
\end{equation}
and the effective divisor $\sum\limits_i \DIV (s_i)$ is normal crossing. 
By the construction of $\tau_k$, we have $ \tau_k \preccurlyeq (\varphi +\psi) $. Therefore 
\begin{equation}\label{add1}
\tau_k \circ  \pi \preccurlyeq (\varphi +\psi)\circ \pi.
\end{equation}
Applying Siu's decomposition of closed positive current theorem to $dd^c (\varphi\circ \pi)$, $dd^c (\psi\circ \pi)$ respectively, \eqref{add1} and \eqref{add2} imply the existence of numbers $a_i,  b_i \geq 0$ satisfying:

$(i)$  $a_i +b_i = c_i$ for every $i$.

$(ii)$ $\sum\limits_i a_i \ln |s_i|  \preccurlyeq   \varphi\circ \pi$ and $\sum\limits_i b_i \ln |s_i| \preccurlyeq   \psi\circ \pi$.

Let $p\in \mathbb{N}$ be an arbitrary integer,  $J$ be the Jacobian of $\pi$, $x\in X$, $f\in \mathcal{I} (p\varphi)$ and $g\in \mathcal{I} (p\psi)$.
$(ii)$ implies that
\begin{equation}\label{integralcondition}
\int_{\pi^{-1}(U_x)} \frac{|f\circ \pi|^2 |J|^2}{\prod\limits_i |s_i|^{2 p a_i }} < +\infty \qquad\text{and}\qquad \int_{\pi^{-1}(U_x)} \frac{|g\circ\pi|^2 |J|^2}{\prod\limits_i |s_i|^{2 p b_i}} < +\infty
\end{equation}
for some small open neighborhood $U_x$ of $x$.
Since $\sum\limits_i \DIV (s_i )$ is normal crossing, \eqref{integralcondition} implies that
$$\sum_i (p a_i -1)\ln |s_i|  \preccurlyeq \ln (|f\circ \pi  |) +\ln |J| \qquad\text{and}\qquad
\sum_i (p b_i -1)\ln |s_i|  \preccurlyeq \ln (|g\circ \pi  |)+\ln |J|.$$
Combining this with $(i)$, we have
\begin{equation}\label{add3}
\sum_i (p c_i -2)\ln |s_i|  \preccurlyeq \ln (|(f\cdot g)\circ \pi  |) +2\ln |J| .
\end{equation}
Note that $J$ is independent of $p$, and $c_i >0$. \eqref{add3} implies thus that,
when $p \rightarrow +\infty$, we can find a sequence $\epsilon_p \rightarrow 0^+$, such that
\begin{equation}\label{singcompare}
\sum_i p c_i (1-\epsilon_p)\ln |s_i|  \preccurlyeq \ln |(f\cdot g)\circ \pi  |  .
\end{equation}

Since $f$ (respectively $g$) is an arbitrary element in $\mathcal{I} (p\varphi)$ (respectively $\mathcal{I} (p\psi)$), by the constructions of $\varphi_p$ and $\psi_p$, 
\eqref{singcompare} implies that
$$ \sum_i c_i  (1-\epsilon_p) \ln |s_i|    \preccurlyeq  (\varphi_p +\psi_p ) \circ \pi.$$
Combining this with the fact that $ (1-\epsilon_p) \tau_k \circ \pi\sim\sum\limits_i c_i  (1-\epsilon_p) \ln |s_i| $, we get 
$$ (1-\epsilon_p) \tau_k \circ \pi \preccurlyeq  (\varphi_p +\psi_p ) \circ \pi .$$
Therefore $ (1-\epsilon_p ) \tau_k   \preccurlyeq  (\varphi_p +\psi_p ) $ and \eqref{equivalproof} is proved.

\end{proof} 

\noindent{\bf Acknowledgement.} I would like to thank Professor S. Boucksom for calling my attention to this problem and for further helpful discussions. I would also like to thank  Professor J.-P. Demailly for helpful remarks and encouragement.

\Addresses

\end{document}